%% file: main.tex
\documentclass[a4paper]{amsart}

\usepackage[foot]{amsaddr}

\usepackage[utf8]{inputenc}
\usepackage[T2A]{fontenc}
\usepackage[russian,english]{babel}

\usepackage{url}
\usepackage[hidelinks]{hyperref}
\usepackage{comment}
\usepackage[left=3.5cm,right=3.5cm,bottom=4cm]{geometry}
\usepackage{needspace}
\usepackage{pgfornament}
\usepackage{enumitem}
\setenumerate[0]{label={(\roman*)}}

\theoremstyle{theorem}

\theoremstyle{definition}
\newtheorem*{definition}{Definition}

\newcommand{\CC}{\mathbb C}
\newcommand{\eps}{\varepsilon}
\newcommand{\TLA}{$\text{TLA}^+$}
\newcommand{\ZZ}{\mathbb Z}
\newcommand{\Lesliedots}{..}
\newcommand{\Lesliedef}{\stackrel\Delta =}
\newcommand{\Yuridoi}[1]{\href{http://dx.doi.org/#1}{DOI: \path{#1}}}

\newcommand{\secauthor}[1]{{\vspace{-.5\baselineskip} \centering \it by #1\par}\vspace{\baselineskip}}

\begin{document}

\title{Mathematical Proof Between Generations}

% AUTHOR COMMANDS
%\author{}
%\address{}
%\curraddr{}
%\email{}
%\thanks{}

\newcommand{\affil}[1]{\textsuperscript{\,(\lowercase{#1})}}

\author[J. Bayer et al.]{Jonas Bayer\affil{a}}
%\address{Freie Universität Berlin}

\author[]{Christoph Benzmüller\affil{b,\,a}}
%\address{Otto-Friedrich-Universität Bamberg and Freie Universität Berlin}

\author[]{Kevin Buzzard\affil{c}}
%\address{Imperial College London}

\author[]{\mbox{Marco David\affil{d}}}
%\address{École Normale Supérieure de Paris}
\address{\vspace{-\baselineskip}}
\thanks{\textit{Correspondence:}  \href{mailto:proof-between-generations@ens.fr}{proof-between-generations@ens.fr}}

\author[]{Leslie Lamport\affil{e}}
%\address{Microsoft Research}

\author[]{Yuri Matiyasevich\affil{f}}
%\address{Steklov Intitute of Mathematics}

\author[]{\mbox{Lawrence Paulson\affil{g}}}
%\address{University of Cambridge}

\author[]{Dierk Schleicher\affil{h}}
%\address{Aix-Marseille Université}

\author[]{Benedikt Stock\affil{i}}
%\address{University of Oxford}

\author[]{Efim Zelmanov\affil{j}}
%\address{University of California, San Diego}

\maketitle

\begin{quote}
\vspace{-\baselineskip}
\footnotesize
\textsc{Affiliations.}
\begin{enumerate}[label={\textsuperscript{(\alph*)}}]
    \item Freie Universität Berlin
    \item Otto-Friedrich-Universität Bamberg
    \item Imperial College London
    \item École Normale Supérieure de Paris
    \item Microsoft Research
    \item Steklov Institute of Mathematics at St. Petersburg
    \item University of Cambridge
    \item Aix-Marseille Université
    \item University of Oxford
    \item University of California, San Diego
\end{enumerate}
\end{quote}

%\vspace{2\baselineskip}
\begin{abstract}
A \emph{proof} is one of the most important concepts of mathematics. However, there is a striking difference between how a proof is defined in theory and how it is used in practice. This puts the unique status of mathematics as exact science into peril. Now may be the time
to reconcile theory and practice, i.e. precision and intuition, through the advent of \emph{computer proof assistants}.
For the most time this has been a topic for experts in specialized communities. However, mathematical proofs have become increasingly sophisticated, stretching the boundaries of what is humanly comprehensible, so that leading mathematicians have asked for formal verification of their proofs. At the same time, major theorems in mathematics have recently been computer-verified by people from outside of these communities, even by beginning students.
This article investigates the gap between the different definitions of a proof and possibilities to build bridges. It is written as a \emph{polemic} or a \emph{collage} by different members of the communities in mathematics and computer science at different stages of their careers, challenging well-known preconceptions and exploring new perspectives.
\end{abstract}

\vspace{2\baselineskip}

\setcounter{tocdepth}{1}
%\tableofcontents

% === MAIN CONTENT ===

% INTRODUCTION OF THE DILEMMA OF MATHEMATICS
\input{1-schleicher}

% HILBERT 10: WHY FORMALIZE ESTABLISHED RESULTS
\input{2-matiyasevich}

% PERSPECTIVE OF RESEARCH MATHEMATICIAN
\input{3-zelmanov}

% MAKING MATH MORE RIGOROUS (WITHOUT FORMALIZING IT ALL)
\input{4-lamport}

% INTRODUCTION TO COMPUTER ASSISTED PROOFS / ACCOMPLISHMENTS
\input{5-benzmueller}

% HILBERT 10: THE CASE STUDY
\input{6-david}

% TEACHING FORMAL SYSTEMS
\input{7-buzzard}

% CONCLUSION, OUTLOOK
\input{8-paulson}

%\begin{acknowledgment}{Acknowledgments}
%\end{acknowledgment}

\vspace{3\baselineskip}

\begin{center}
\pgfornament[width=.7\linewidth,color=gray]{88}
\vspace{3\baselineskip}
\end{center}

\bibliographystyle{abbrv}
\bibliography{common-references,1-references,2-references,3-references,4-references,5-references,6-references,8-string,8-atp,8-funprog,8-general,8-isabelle,8-theory,8-crossref}

\end{document}

%% file: 1-schleicher.tex
\needspace{4\baselineskip}
\section{What is a mathematical proof?}\label{sec:schleicher}
\secauthor{Dierk Schleicher}

Mathematics often prides itself as the most fundamental of all
sciences: a mathematical truth, once established, will be true
forever. But what exactly is a mathematical proof? Here is one
possible answer in terms of the other fundamental concept in
mathematics, a definition:

\begin{definition}[Proof: formal definition]
  A mathematical proof is a sequence of arguments that are based on a
  given set of axioms and formally deduce consequences following
      formal rules of deduction.
\end{definition}

One might know this definition of proof from high school proofs, where
one shows that the sum of the angles in any Euclidean triangle is
equal to a straight angle, or perhaps gives a formal proof of the
Pythagoras theorem. But sooner rather than later, one discovers that
this definition is too clumsy and impractical for any profound result
in mathematics. In practice, mathematicians use a very different
definition of proof.

\begin{definition}[Proof: practical definition]
  A mathematical proof is a sequence of arguments that convinces an
  educated reader.
\end{definition}

So one might object that mathematicians routinely make one of the most
fundamental mistakes: define things in one way and then use them in a
different way. But does this not give away all the fundamental virtues
of mathematics?

Many mathematicians are quick to point out that a typical proof, given
as a convincing sequence of arguments, can be elaborated in more
detail: in principle, for every step one may ask why this is so, and
one can insert additional steps of explanation. Each of these steps,
in turn, may be expanded if need be, until eventually one arrives at
the axioms. One would expect that any of these expansion steps can be
performed by the author of a proof, until one arrives at sufficiently
fundamental levels when one reaches statements that are already
established. It is a common understanding among mathematicians that
most current mathematical knowledge can be thus justified -- in
principle.

But it is true that mathematical ``proofs'' written by humans
occasionally do have errors. Many of them are of a trivial kind when
some missing steps can be inserted, or an overlooked special case can
be treated in the same way as the rest of the proof. Sooner or later,
one is often told, any mistakes would be found by the scrutiny of the
mathematical community.

But now and again there are theorems that were thought to be proven
where much later the proof is recognized to be false or incomplete,
and where perhaps even the result itself is recognized as wrong. How
would you know that for a given theorem the ``sooner or later'' might
not happen tomorrow, when a possibly fundamental flaw is discovered?
How can one then be certain of the mathematical correctness of a
proof?

A rather prominent classical example is the Four Color Theorem: this theorem was thought to be proved already in the late 19th century, when not one but two proofs were given --- that were both found to be incorrect eleven years later. 

\subsection*{A recent example from dynamical systems.}

Proofs that are later discovered to be flawed are not just a theoretical possibility, or an anecdote from the early times of mathematics (when perhaps standards of mathematical rigor were not so high). This is illustrated in a very current example from
dynamical systems. The claim is stated quite easily:

\emph{%
  Let $f\colon\CC\to\CC$ be a holomorphic function but not a
  polynomial (i.e.\ a transcendental entire function). Then $f$ can
  have at most one completely invariant domain; that is a domain $U$
  with $f(U)=U=f^{-1}(U)$.%
}
Until very recently, this result would be stated as a theorem, proved by Noel Baker, an eminent pioneer in the field, in 1970~\cite{Baker70}. However, Rempe-Gillen~\cite{Lasse_Dave} observed about a half century later that the proof is flawed and cannot be fixed by the methods stated. In particular, a key step in the proof is the statement that if
$U,V\subset\CC$ are disjoint simply connected domains, then one of
$f^{-1}(U)$ and $f^{-1}(V)$ must be disconnected. This statement,
however, is false; even as simple a function as $f(z)=e^z+z$ provides
a counterexample: it has infinitely many disjoint simply connected
domains with connected preimages. Baker's proof thus cannot be
repaired along the lines of its original version. The main result, as
of today, is an open question --- once again, half a century later.

One might wonder how relevant this result is --- perhaps the error went
unnoticed because nobody cared? Unfortunately, this is not so. Rempe-Gillen \cite[Section~9]{Lasse_Dave} gives a shocking list of results that
depend on the flawed paper in a variety of ways. One list contains
several publications by a variety of authors that all depend either on
the flawed paper, or even on other papers that used the flawed result,
but where the main results can be fixed by other methods developed in
\cite{Lasse_Dave}. Another list contains several publications, some 35
years old, that use the flawed paper and for which the main results
must now be considered open once again. Yet another list contains two
survey papers, between 25 and 30 years old, that refer to the result
in \cite{Baker70}. And perhaps most embarrassingly, there are several
publications, some of which meanwhile classical and much-cited, where
the flawed method of proof was found so useful that it was adapted,
further developed, and generalized (without recognizing the flaw). As
a result, a whole area of mathematics has had to regroup and sort out
how its theory was affected by a rather ``convincing'' mistake made
half a century earlier.

This is just a recent example, rather embarrassing for the
mathematical community, that shows that problems can be discovered at
any time even in supposedly well-known results. 

However, it may well be that mathematics is now developing
towards the point when, finally, the two definitions of ``proof''
can be reconciled. How this may come about is one of the key topics of this paper.

Another development in current mathematics, leading to similar conclusions, is that proofs tend to become longer and more complex. For instance, Peter Scholze has recently observed that proofs have become so complex that some could no longer be stored as a whole in human memory, nor be verified by referees. He thus asked that one of his key results on condensed mathematics be formally verified, a project which has attracted broad attention~\cite{liquid_tensor_nature, hartnett2021}.
A number of years earlier, Vladimir Voevodsky had similarly been worried about the correctness of results in his field leading to the development of ``Homotopy Type Theory''~\cite{voevodsky-origins}.

\medskip

One of the inspirations for the paper was a workshop at the \emph{Heidelberg Laureate Forum} in 2018, titled ``The future of mathematical proof'', where many of the coauthors of this paper were present, and where some of the aspects treated here were discussed and developed.

%\textcolor{red}{Other leading mathematicians have started questioning the status quo over the past two decades. Prominently, Vladimir Voevodsky's program on univalent foundations \cite{voevodsky-origins} was prompted by an ever growing complexity of mathematical proofs. Recently, the call by Peter Scholze to formally verify one of his proofs has again attracted broad interest \cite{liquid_tensor_nature, hartnett2021} in this debate.}

This paper should be seen as a kaleidoscope of many aspects of mathematical proofs and computer proof assistants, contributed by various researchers with diverse backgrounds. In Section~\ref{sec:matiyasevich}, \emph{Yuri Matyiasevich} describes his vision that, soon, mathematical publications need to be supported by formal proofs, and that this may well be accomplished by the coming generation. In Section~\ref{sec:zelmanov}, \emph{Efim Zelmanov} argues in favor of the meaning and explanations of mathematics, rather than formal verification. In Section~\ref{sec:lamport}, \emph{Leslie Lamport} introduces a way of presenting mathematical proofs in a way that ``should make it much harder to publish false proofs''. In Section~\ref{sec:benzmueller}, \emph{Christoph Benzm\"uller} discusses several major accomplishments of proof assistants and derives a vision for integrated formal and traditional proofs. Then, in Section~\ref{sec:david}, \emph{Jonas Bayer, Marco David, and Benedikt Stock}, three students and thus members of the coming generation, describe how they learned to work with proof assistants from early on, and completed perhaps the first major formalization done outside of the expert community. In Section~\ref{sec:buzzard}, \emph{Kevin Buzzard} explains how he makes proof assistants ``sexy'' to mathematicians and recounts how he incorporated them into classroom teaching. In Section~\ref{sec:paulson}, \emph{Lawrence Paulson} outlines open challenges and a perspective for the future.

%% file: 2-matiyasevich.tex
\section{Why formalize mathematical results? And why Hilbert 10?}\label{sec:matiyasevich}
\secauthor{Yuri Matiyasevich}

A decade ago, the St.\,Petersburg mathematical society
 held a meeting titled
``Mathematical proof: yesterday, today, tomorrow''
\cite{matobvo1,matobvo2,matobvo3}. Being one of the
three spokesmen, I completely disagreed with the
previous speaker \cite{matobvo1} (see also \cite{Vavilov}),
and dared \cite{matobvo2,compinstr} to  publicly
announce
the following:

\begin{center}
\vspace{\baselineskip}
\fbox{\parbox{120mm}{
  \begin{center}
  PREDICTION
  \end{center}

In 25 years mathematical journals (if they are destined to survive
so long) won't take into consideration any paper unless
it is accompanied by proofs which can be verified
by a computer.

}}
\vspace{\baselineskip}
\end{center}

Already at that time at least one (and only one?)
journal existed
in which proofs passed preliminary computer verification.
This was \emph{Formalized Mathematics}
\cite{fm} founded in 1990.
Over the past 10 years, mathematical journals
continued to multiply in number but I do not know
of any new journal requiring formalized proofs.
Nevertheless, I dare to repeat my prediction
\emph{verbatim et litteratim}, that is meaning
``In 25 years \emph{from now} ...''.

In fact, the progress in computer verification
of proofs is very impressive, both in the software development
and in the number and the significance of
actually verified theorems. You can consult \cite{100} for
a compendium of such achievements; the list of proofs is
regularly updated  but initially restricted
by a selection of ``Top 100''  theorems.

But what could be the goal of computer verification of
results already proved by
human-beings? One obvious answer is as follows: to get
additional trust in the correctness of a proof
(there are many examples of important and widely used
theorems with flaws in their proofs
which remained unrevealed for decades). However, some people believe
that computers themselves
are not sufficiently reliable because
of possible errors in their design/software/runs.
  In my opinion ``in 25 years'' the progress
  in this area will remove
  such objections completely.

Formalized proofs are vital
for the very ambitious project
of a \emph{World digital mathematics library}
 \cite{digmath}. The ultimate goal of this project is
 to transfer all \emph{mathematical knowledge}
 (not just axioms/definitions/theorems/proofs)
 to computers. This would require tremendous efforts
 but who would care to spend her/his precious
 time on the meticulous presentation of known results?
 Mathematicians prefer to produce new
 ideas (definitions, hypotheses, theorems)
 and do not   appreciate  the hard work
 of writing proofs with all minute details
 (luckily, computer scientists do appreciate
 such occupation).

 The following remedy was proposed a long time ago:
 senior, ``retired'' mathematicians, who are no longer
 capable of generating new brilliant mathematical
 ideas, could devote the rest of their lives to
 ``teaching mathematics to computers''.
 But there is an opposed option: this could be done
 by the young people who  just start to  master
 mathematics. In this way they can get an acute feeling
 of what mathematical rigor is.
 However, the following remains uninvestigated:
 how would involvement in such an activity influence
 the ability to create new mathematics?

  So when I heard that a group of students from Jacobs
  University in Bremen had studied, under supervision
  of Professor Dierk
  Schleicher, the proof of the undecidability of
  Hilbert's tenth problem, I suggested to them to
  demonstrate their understanding of the whole construction
  by producing a fully formalized proof.
  My role in the project was very restricted: I
  supplied the students with a sufficiently (for human-beings)
  detailed proof, and I was responsible for the
  choice of
  \emph{Isabelle} as the verifier. I was happy to see
  that the students  became very enthusiastic,
  maybe because they (and me too)
  at that moment underestimated the amount of required  work.

Hilbert's tenth problem is not very difficult for formalization.
It is a bit strange that at first we were  waiting for it
for half a century but then
four independent verifications \cite{dprm_isabelle,Lean10,form3,form1a}
emerged in a short time. 

%% file: 3-zelmanov.tex
\section{On proof and progress in mathematics: the perspective of a research mathematician}\label{sec:zelmanov}
\secauthor{Efim Zelmamov}

I will add my 5 cents to the wonderful discussion of computer proof verification, organized by students. For more than 40 years, I’ve lived in the world of proofs and, sometimes, complicated proofs.

If I were told that a proof is correct because a computer program says so, but I don't see big ideas "turning the wheels", then probably I would continue thinking about the problem, as if a computer blessing did not exist.

The purpose of a proof is \emph{understanding}. For mathematicians it is not enough to know if this or that statement is correct or not. They want to know \emph{why} it is correct or not. Often this understanding comes as a link to some big ideas that come into play now and then in different contexts. I am not able to say it better than W. Thurston in his beautiful paper on proofs (see \cite{Thurston94}).

In my opinion, a credible computer verification of a proof is an amazing achievement in AI. It is valuable and interesting for its own sake, leaving alone proofs. It may also find other applications.

I have to admit that as far as a straightforward computation is concerned, I trust computers more than humans.

Difficult and important proofs are often written at the edge of human intellectual abilities (think of the proof of G. Perelman of the Poincaré Conjecture). Should we expect the authors to present it in a computer-friendly form? I am afraid that it is too much to ask.

Modern mathematics created a new challenge: proofs of enormous complexity. I know proofs of some far-reaching statements that have been around for quite a while and still nobody (except the authors) can say that they understand all the details. The only hope is that good proofs are like live organisms. They evolve with time and go through a natural selection. New ideas appear and bring new understanding. And new proofs may appear.

I will end with the controversial statement that in practice a proof is what is considered to be a proof by all mathematicians.

%% file: 4-lamport.tex
% FROM THE PREAMBLE
% \newcommand{\eps}{\varepsilon}
% \newcommand{\TLA}{$\text{TLA}^+$}
% \newcommand{\ZZ}{\mathbb Z}
% \newcommand{\Lesliedots}{..}
% \newcommand{\Lesliedef}{\stackrel\Delta =}

\section{Making math more rigorous}\label{sec:lamport}
\secauthor{Leslie Lamport}

Mathematics, as practiced by most mathematicians, is not very rigorous. There is evidence that about $1/3$ of all published, refereed math papers contain significant errors --- incorrect proofs or theorems that their authors believed to be correct. (I have presented evidence elsewhere \cite{Lamport_ErrorsInProofs}.) Math can be made more rigorous, and mathematicians can make fewer errors, by replacing archaic customs with more sensible practices. Here is how.

\subsection*{Formulas.} 
A few hundred years ago, formulas were written in prose. Today, mathematicians recognize the advantages of writing formulas in mathematical notation: they're shorter, easier to understand, and easier to manipulate. Replacing prose by mathematics must have reduced errors. 

Mathematicians think they've stopped using prose to write formulas. 
They're wrong. They've replaced only some of the prose in their formulas 
by math. Consider this definition of what it means for $\lim_{x\to a} f (x)$ to equal b. 
\begin{align}
&\text{For all $\eps>0$ there exists $\delta >0$ such that, for all $y$, if $0<|y-x|<\delta$} \nonumber \\ 
&\text{then $|b - f (y)| <\eps$. }
%\nonumber
\label{Eq:LeslieLimit}
\end{align}
A mathematician would find \eqref{Eq:LeslieLimit} perfectly normal, even though it's a mathematical formula written with many words. Here is that formula written without words:
\begin{equation}
\forall \eps>0 \colon \exists\delta>0 \colon \forall y \colon (0<|y-x|<\delta) \implies (|b-f(y)|<\eps) 
\label{Eq:LeslieLimit2}
\end{equation}
I believe most mathematicians would find \eqref{Eq:LeslieLimit2} harder to understand and uglier than \eqref{Eq:LeslieLimit}. I expect mathematicians a few hundred years ago would have found $0 < |y - x | < \delta$ hard to understand and ugly. 

Why write \eqref{Eq:LeslieLimit2} rather than \eqref{Eq:LeslieLimit}? For the same reason we don't write \emph{0 is 
less than the absolute value of\dots}: It's shorter, easier to understand (when 
you become comfortable with the notation), and easier to manipulate. And 
it will reduce errors. Show elementary calculus students the definition \eqref{Eq:LeslieLimit} 
and ask them to write what it means for $\lim_{x\to a}f(x)$ equals $b$ to be false.  
I doubt if many of them would get it right. Teach them a little elementary logic and they could easily compute the negation of \eqref{Eq:LeslieLimit2}. The most obvious use of words in formulas is to express logical operations; but they are also used in other ways, such as describing sets and functions. 

Formulas written without words can now be manipulated by computer programs. Programs can easily compute the negation of \eqref{Eq:LeslieLimit2}. They can't compute the negation of \eqref{Eq:LeslieLimit}.\footnote{Restricted, unnatural languages have been proposed for writing formulas approximately like \eqref{Eq:LeslieLimit} so they can be understood by a computer program. Such languages are of little or no use to people not afraid of mathematics. 
} Those programs can help students become comfortable with mathematical concepts, if the concepts are described with math rather than prose. 

Mathematicians think it's difficult to write formulas completely mathematically, without words. I have asked a number of mathematicians how long a completely rigorous, wordless definition of the Riemann integral would be --- assuming definitions of the set of real numbers and its arithmetic operations, as well as simple set theory. I've received answers ranging from 50 lines to 50 pages. They're wrong. 

I've developed a language called \TLA\ that engineers use to write completely formal mathematical descriptions of computer systems. It has tools for checking the correctness of their mathematics. The Riemann integral can be defined in \TLA\ in about a dozen lines. 

\subsection*{Proofs.}

A few hundred years ago, proofs were written in prose. They still are. Mathematicians haven't even begun to change the way they write proofs. They think their proofs express rigorous logical reasoning. They're wrong. Their prose proofs are written in a literary style that obscures the logic of the proof. Consider the following opening sentence of a proof from an elementary calculus book by Michael Spivak \cite[page~170]{Spivak} --- a book that is considered to be very rigorous. 
\[
\text{Let $a$ and $b$ be two points in the interval with $a < b$.} 
\]
It is obviously wrong because the interval in question could consist of a single point, so it might be impossible to choose $a$ and $b$. That sentence is actually part of a correct proof, but the reader must discover for herself the proof hidden inside Spivak's prose.
 
Writing proofs with prose leads to errors. How can those errors be avoided? Most mathematicians and computer scientists believe that the only way is to write machine-checked proofs. This requires writing formulas in a formal language. \TLA\ is simple enough that mathematically unsophisticated engineers can use it, and it is enough like everyday math that mathematicians should find it fairly natural. But it's too simple to be adequate 
for writing the kinds of proofs found in most math journal articles. Formalizing such proofs requires a language too complicated for most engineers to learn --- one that I believe most mathematicians would find quite obscure. Few mathematicians would go to the effort of learning such a language unless it made writing their proofs significantly easier. Today, it makes writing most proofs much more difficult. Routine machine-checked proofs are now practical in just a few situations, including some safety-critical applications. I don't expect this to change in the next couple of decades. 

Fortunately, there is a simple method that anyone can use now to write proofs with fewer errors. It can't eliminate all errors, but it can make them much less likely to occur. Its basic idea is to replace the linear order of ordinary prose by a hierarchical structure, and to name hypotheses and proved facts so they can be referred to later in the proof. Here is a brief explanation of the method; a complete description has appeared elsewhere \cite{Lamport_HowToWriteProof,Lamport_HowToWrite21}. 

A theorem consists of a statement together with its proof. A proof is either a short paragraph or a sequence of statements and their proofs. At each point in a proof, there is a current goal and a set of usable facts that can be assumed in proving that goal. Statements can be written in prose or in math. When written in math, the logical structure of the statement often determines the hierarchical decomposition of its proof. Figure~\ref{Fig:Leslie1} shows the structure of part of a proof containing the statement $A\wedge B \implies C$, in which $C$ is proved by first proving statements $D$ and $E$. Usually, those two statements would easily imply $C$, making the proof of {\sc qed} step $\langle 3\rangle 4$ simple. The number $\langle2\rangle 3 $ indicates that it is the third statement in the level-2 proof of a level-1 statement. 

\newcommand{\customsize}[0]{\footnotesize} %{\footnotesize\sf}
\begin{figure}[htbp]
{\raggedright
\rule{0pt}{1pt}\quad $\langle 2\rangle 3$. \ $A\wedge B\implies C$\\
{\customsize Current goal set to $A\wedge B\implies C$}\\
\rule{0pt}{1pt}\quad\quad$\langle 3\rangle$1. \ \textsc{suffices assume \ $A$, $B$ }\\
\rule{0pt}{1pt}\phantom{\quad\quad$\langle 3\rangle$1. \ \textsc{suffices}} \textsc{prove} \ \ \ $C$\\
\rule{0pt}{1pt}\quad\quad\quad\, Proof: By simple logic. {\customsize Trivial proof that assuming $A$ and $B$, then proving $C$,}\\
\rule{0pt}{1pt}\rule{2pt}{0pt}\phantom{\quad\quad\quad\, Proof: By simple logic.} {\customsize proves the current goal} \\
{\customsize Current goal set to $C$; and $A$ and $B$ added to usable facts.}\\
\rule{0pt}{1pt}\quad\quad$\langle 3\rangle$2. \ $D$\\
\rule{0pt}{1pt}\quad\quad\quad\, Proof of $D$\\
{\customsize $D$ added to usable facts.} \\
\rule{0pt}{1pt}\quad\quad$\langle 3\rangle$3. \ $E$\\
\rule{0pt}{1pt}\quad\quad\quad\, Proof of $E$\\
{\customsize $E$ added to usable facts.} \\
\rule{0pt}{1pt}\quad\quad$\langle 3\rangle$4. \ QED\\
\rule{0pt}{1pt}\quad\quad\quad\, Proof of $C$\\
{\customsize Current goal and usable facts same as before $\langle 2\rangle 3$ except with fact $A\wedge B\implies C$ added.}\\ 
\rule{0pt}{1pt}\quad $\langle 2\rangle 4$. \dots\\
}

% The caption did not want to center itself for some reason...
\caption{A statement and its structured proof.\qquad\phantom{abcabcabc}}\label{Fig:Leslie1}
\end{figure}
This is a straightforward proof, and presented in this way there seems no reason to structure it. But suppose it were a small part of a large proof, and the proofs of $D$ and $E$ were each half a page long. If the proof were written as prose, how could the reader keep track of where the scope of the hypotheses $A$ and $B$ ended, and where it was no longer valid to use $D$? Mathematicians try to handle complexity by using Lemmas; but that just provides one level of hierarchy, which doesn't get you very far. 

Making a proof more rigorous requires filling in all the gaps that could conceal errors. This means making it longer. Making a prose proof longer makes it harder to read. But with hierarchical structure, the extra length makes the proof easier to read. The additional explanation appears at lower levels of the hierarchy, so it doesn't obscure the structure of the proof. This will be especially true when mathematicians stop producing pictures of print on dead trees and start using hypertext, so lower-levels of the proof can be hidden when not being read. 
Avoiding errors requires more detailed proofs than are currently found in journals. Until journals use hypertext, this means writing a detailed proof
to catch errors, then shortening it for publication. That's easy to do with structured proofs: you just replace the lower levels of the hierarchy with short proof sketches. (One can even write \LaTeX\ macros so a single file can produce either version by changing a few characters.) 

Students can learn to write structured proofs by teaching them to write very simple machine-checked proofs in some field. Any good proof system should allow hierarchical structuring of proofs. The language for writing theorems should be simple --- not the kind of complicated language needed for serious math. \TLA\ and its proof system are not ideal, but they could be used if nothing better is available. 

Students should understand that the facts they learn in their math classes can, in principle, be formally proved from simple axioms and proof rules. In practice, we only carry proofs down to the level where we believe the reader will find the steps to be obviously true. That level rises with education and experience. We also sometimes take shortcuts by writing formulas with words. But students and mathematicians should have the confidence that they could make their proofs completely rigorous and carry them down as close as they want to basic axioms. 

Writing hierarchically structured proofs can help you avoid errors; it can't guarantee that you will. You have to be honest with yourself about 
what's obvious and what should be proved. My advice is to write the proof down to a level at which you think everything is obvious, and then go one level deeper. But if you don't care whether your proofs are correct, nothing short of having to write a machine-checked proof will keep you from making errors. 

\subsection*{What Should You Do Now?}

If you agree that writing formulas with words or that writing proofs with prose is silly, just stop doing it. You don't have to wait for others to change. 
\subsubsection*{Formulas.} 
You needn't remove all the words from your formulas. Start by using the quantifiers $\forall$ and $\exists$. Then try eliminating ``$\dots$'', which is not a mathematical operator. The sequence $x_1, \dots, x_n$ is just a function $x$ with domain $1\Lesliedots n$ that maps each $i$ in its domain to $x_i$.\footnote{``\Lesliedots'' is a mathematical operator, defined by $i\Lesliedots j \Lesliedef \{k\in\ZZ\colon i\le k\le j\}$, where $\ZZ$ is the set of all integers. 
} Often, though not always, the math becomes simpler and more elegant if you eliminate the ``\dots'' and instead use the function $x$. Give it a try. Be aware of when you're using words and sloppy notation instead of being rigorous. If you're open to change, you will find that the mathematically rigorous approach is often the simplest. 
If you're a teacher, your students should have learned, or should be learning, the basic math needed to write formulas with fewer words than they now use. Help them to become more comfortable with proper mathematical notation by using it in your classes. 

\subsubsection*{Proofs.} 
There is no reason not to start writing structured proofs now. It takes only a sentence or two to explain to readers how to read them. I've been doing it for about 30 years, and no editor or referee has complained about my proofs. Start by reading how I write structured proofs, but feel free to modify my style as you see fit. There are just two features that should be preserved: hierarchical structuring and the ability to name and refer to hypotheses and already proved statements. 

If you're a professor, teach your students to structure their proofs the way you do. They're not yet set in their ways, and they'll appreciate how 
the structure makes your proofs easier to understand. Encourage them to write structured proofs in all their courses. Other professors are unlikely to complain that the proofs are too rigorous; and they might even be inspired to write them themselves.

%% file: 5-benzmueller.tex
\needspace{4\baselineskip}
\section{What is a proof? What should it be?}\label{sec:benzmueller}
\secauthor{Christoph Benzm{\"u}ller}

Does the notion of a mathematical proof refer to the rigorous but typically rather non-intuitive \emph{formal derivation of a new ``truth'' from its premises using accurately defined rules of inference}? Or is it an \emph{artful communication act} in which the beautiful structures underlying a new mathematical insight are revealed to peers in such a way that they can easily \emph{see} and accept it, and even gain further inspiration?

The former notion of a \emph{formal proof} is primarily concerned with logical rigor and soundness. Intuition and beauty is a secondary concern, if at all. Formal proofs have recently attained increased, albeit quite controversial, attention in mathematics, triggered e.g.~by successful applications of modern theorem proving technology to challenging mathematical verification and reasoning tasks. Some of the settled problems are of such a kind that human cognition alone has apparently reached its limits for attacking them. Respective examples include:

\renewcommand{\theenumi}{(\roman{enumi})}
\begin{enumerate}%[leftmargin=*,labelindent=0pt]
 	\item \emph{The Four-color theorem:} this notorious challenge was solved already in 1977 using automated theorem technology by Appel \& Haken  \cite{1977SciAm.237d.108A}, and an interactive formal proof was more recently developed by Gonthier \cite{Gonthier08} in the proof assistant Coq~\cite{Coq}.
    \item \label{hales} \emph{The Kepler conjecture} (about best sphere packings in Euclidean 3-space): a board of experts had surrendered the complex task of reviewing Hales' solution for the Annals of Mathematics. Eventually, Hales and his team mastered it in interaction with HOL Light. As the main result,  a formal proof is now available that is independently verifiable --- by humans and computer programs \cite{hales-formal-Kepler}.
    \item \label{heule} \emph{The Pythagorean Triple Problem} (whether all positive numbers can be colored by two colors so that there is no monochromatic 
    Pythagorean triple): this challenge problem was solved by Heule \& Kullmann \cite{Heule2017} using automated SAT solving technology (boolean \textsc{SAT}isfiability), and the formal proof that was generated by the computer program is of enormous size (about 200TB); it is still  independently verifiable though (at least by machines). This line of research has been continued by an automatic solution for Schur Number Five \cite{Heule2018}, and even Keller's conjecture \cite{Heule2018} was recently resolved.
    \item My own current work with colleagues focuses on higher-order meta-logical reasoning technologies \cite{J41} that have enabled us to detect and explain errors and problems in peer-reviewed publications in mathematics, computational metaphysics, and machine ethics. These include the discovery of an unnoticed inconsistency in G{\"o}del's modern variant of the ontological argument for the existence of God~\cite{C55}, the discovery and clarification of a deeply rooted paradox in Zalta's Principia Logico-Metaphysica~\cite{J50}, and the revelation of some minor problems in a well-known textbook on category theory~\cite{J40}. 
\end{enumerate}

It is to be noted that the formal reconstruction of mathematical work is generally a very resource- and time-consuming task, for example due to the lack of large and easily reusable libraries of formalized mathematics, or software that does not pose additional challenges along the way. Conversely, some technical mathematical results, such as those mentioned in (iii), may not support intuitive and insightful mathematical proofs.
In general, mathematics is confronted with increasingly complex problems whose solution and subsequent evaluation (e.g., by peer review) require techniques that go beyond traditional practice. Examples (i)--(iv) above are only early evidence of this kind. When technologies are available that can help detect errors in publications, they should certainly be used to optimize scientific quality. Formal proofs should therefore take a more central role, in maths and beyond.

In fact, I see it as a societal duty to take up this challenge. Clinging to the traditional notion of mathematical proof alone is not an option in an increasingly technological world. Think of areas such as ``program verification in computer science'' or ``trusted AI'', where ideally we want formal guarantees that implemented, complex solutions are mathematically correct, but where intuitive, traditional proofs might be lacking.

Nonetheless, formal proofs alone are of limited interest and they should ideally always be coupled to intuitive proofs. Explainability, transparency, and intuition must remain virtues of the highest priority, not only in mathematics. In the long run, the increased trustworthiness and beauty of a combined approach, where both notions are coupled pari passu, will justify the additional resources that must now be invested. Publishing errors (even minor ones) will be prevented by formal proofs, and designing ill-conceived and inaccessible theories will be impeded by the demand for mathematical intuition and beauty.

\subsection*{So, what should a proof be?} In conclusion, it should ideally be both a human-oriented traditional proof and a machine-oriented formal proof. Traditional mathematical proofs are made by, and consumed by, humans, while formal proofs are predominantly generated with, and consumed by, machines. Yet, by pairing both, the antipodes will increasingly engage in a harmonious dance.

One step further, beyond \emph{pairings}, one may dream of the \emph{integration} of mathematical and formal proofs into one object. Modern proof assistants such as Isabelle/HOL provide increasingly intuitive languages for constructing and representing proofs, but significant scientific progress is still needed to achieve this ambitious goal. Finding a comprehensive solution that integrates both notions of proof can even be understood as a challenge for AI, as it requires a seamless semantic integration of natural language, diagrams, formula language, etc., up to the exchange of meaningful arguments between humans and theorem provers.

Recent work by Marco David, Benedikt Stock, Jonas Bayer and their fellow students gives reason for hope. Their verification project \cite{Digit_Expansions-AFP, dprm_isabelle, DPRM_Theorem-AFP} of Hilbert's $10^{\operatorname{th}}$ problem is significantly smaller in scope than some projects mentioned above, but differs in that the mathematics students began their formalization project, encouraged by Matiyasevich, with no prior knowledge of proof assistant technology. Nevertheless, they rose to the challenge and made great progress by working independently with the proof assistant. As far as I know, this is the first extensive formalization project to come from outside the community, so it is an excellent demonstration of the maturity that modern proof assistant technology has now reached.

%% file: 6-david.tex
\section{Proof assistants: The right way to learn mathematics?!}\label{sec:david}
\secauthor{Jonas Bayer, Marco David, and Benedikt Stock}

Little did we know at the time about Yuri’s ``grand plan'' to establish his vision on the future of mathematical proof when he visited our university back in 2017 and presented the idea to formalize his DPRM-theorem. Little did we know about our role as “guinea pigs” who, in a carefully set-up experiment, should establish the idea that mathematical proof verification by computers is feasible even for young and inexperienced university students --- a role we are now very grateful to have been given. The following paragraphs elaborate what we have learnt in the years since then, now that we are the ones who are supervising another project by a fresh set of students that had never worked on computer-verified proofs before.
 
Yuri set up the scene swiftly; he presented his theorem to us and invited us to work on its computer verification. Neither he nor us had any idea how much effort this would prove to be. In fact, at the time, it seemed to us that, once the proof of the DPRM theorem was understood, we merely needed to ``translate'' its arguments in such a way that Isabelle could verify them. Anyone who has ever worked with an interactive theorem prover, however, knows that the word ``translate'' does not do the process much justice. In reality, this includes filling the possible gaps that are frequently found in mathematical papers. Common phrases such as ``it is easy to see that...'' needed to be brought to logical life during this process. Soon, we realized that the challenges often lay in lemmas that looked rather innocent. So much so that we were forced to five times reconsider our formalization of the mere concept of a register machine before our formal definition proved useful.

Despite the challenges we encountered, we insist that learning how to work with an interactive theorem prover is fully feasible, although not yet easy. In a previous paper \cite{bayer-beginners-quest}, we have reflected upon our own mistakes and the challenges we encountered. A key realization is that making progress is almost impossible without an expert at hand to answer questions. Hence, to popularize proof assistants among the next generation, we are now proactive ourselves in mentoring a new group of students learning how to formalize mathematics.
The current need to pass on experience person-to-person distinguishes proof assistants from most programming languages or computer algebra systems like Mathematica. Here, online Q\&A forums such as Stack Overflow\footnote{\url{https://stackoverflow.com}} provide a freely accessible and searchable database of almost any imaginable question about these tools, and hence provide the means to make them accessible to the broad public. Interactive theorem provers and their currently secluded communities must follow suit to become useful for the average mathematician!

Beyond the technical skills and results, this project prompted a paradigm shift in the way we now view ordinary mathematics: In our university courses, we no longer ask if a proof is convincing to us, but instead wonder if we would be able to formalize it in Isabelle. As the computer often exhibits assumptions or edge cases that humans gloss over, this thought pattern leads to a more rigorous approach to the (informal) argument. Our new, sharpened perspective was largely forged through the interaction with the computer and its unique manner of reasoning. It illustrates how to reconcile the two conflicting definitions of ``proof'' in the next generation of mathematicians.

%% file: 7-buzzard.tex
\section{Some thoughts on the formalization of mathematics}\label{sec:buzzard}
\secauthor{Kevin Buzzard}

\vspace{-\baselineskip}
\subsection*{Making formalization sexy}
Mathematics has fashions. The Langlands Philosophy seems to have been a fashionable area for many decades now. Conjectures get proved, new conjectures get made, there was the classical theory, and then a mod~$p$ theory (crucial for Fermat's Last Theorem), and now a bewildering $p$-adic theory as well as geometric and categorified versions. I will unashamedly confess that my work with Johan Commelin and Patrick Massot where we translated the line ``Let $X$ be a perfectoid space'' into Lean's language was an attempt to \emph{market} theorem proving software to mathematicians. We wanted to show them that theorem provers are now ready to handle fashionable modern mathematics. No, we cannot do anything spectacular yet like come up with an incomprehensible billion-line long proof of the Riemann Hypothesis (indeed, computers generating proofs of hard mathematical conjectures of mainstream interest is science fiction right now, and may well remain so for decades to come). However, let's think about this. If we stick with the status quo (i.e., essentially nobody types sexy mathematics into proof assistants), then proof assistants will \emph{never} be able to do sexy mathematics, because they will simply never learn it! Computers certainly can't directly read the crap we write in our papers (and humans often can't read them either). Humans need to do the translation by hand, and the sooner the better. So, whose \emph{job} is it to type in the \emph{statements} of the global Langlands conjectures into these systems? Who will do this, and give humanity the chance to make computer teaching and proof mining resources for students learning and working in the Langlands Program? Surely it is our job as mathematicians. Nobody else is going to do it -- we cannot expect computer scientists to become technical experts in the Langlands Philosophy; it is much easier for mathematicians to learn a programming language. If the maths staff won't do the translation, let's get the maths PhD students doing it. If you've recently graduated with a PhD from a pure mathematics department, can you, or even can humanity, \emph{state} in Lean the main theorems you proved in your thesis? The practical answers with current technology are still highly dependent on the nature of the area of the thesis. For my MSc and PhD students I am 100\% certain that we can not only state the main results in Lean, but also prove them.

Looking to the future, what if mathematicians start flocking to formalization, and all of a sudden we have got perfectoid spaces and the statements of the Langlands conjectures in Lean and also in Arend (a HoTT prover) and Isabelle/HOL (a simple type theory prover) and MetaMath (a set theory prover) and Coq and HOL 4 and HOL Light and Mizar and cubical Agda and\ldots. What then? I have already said that, in my opinion, it is science fiction to expect that these systems will start \emph{proving} the Langlands philosophy.

But the following is near-reality. Software such as Lean can be used to power an interactive resource for PhD students learning algebraic geometry or some other dense area. Once we have a database of \emph{statements} of many theorems corresponding to tags in the Stacks Project or Kerodon (online databases of algebraic geometry and other sexy mathematics), we can let the computer scientists take over. They have tools known as hammers, which can attempt to build proofs or counterexamples to propositions fed to the system, using the database we mathematicians constructed. Note that such a tool \emph{does not need any formalized proofs} and so even a huge database like the the statements of the theorems in the Stacks project could be constructed manually over a period of several person-years, ideally by young algebraic geometers interested in trying new ways of learning the area.

\subsection*{Making formalization fun}

The examples in the previous paragraph show that the mathematical community might well come to benefit from having serious mathematics formalized in a theorem prover. Experience suggests that it is mostly young people who formalize. Thus teaching young people about formalization is important; however it is a different task to teaching young people about mathematics. The Natural Number Game is a browser-based Lean game that came out of many hours at Imperial just writing random undergraduate-level mathematics puzzles in Lean and then watching students solving them. The idea of building up facts about natural numbers from first principles was a big hit, and it went on from there. Students say to me ``I've finished the natural number game, what next?'' The correct answer to that is ``Install Lean following the instructions on the community website''. But after they've done that, what next? This depends on their mathematical interests and abilities. I often encourage an enthusiastic student to formalize some mathematics they already know well for their first Lean project. The rule is: if it compiles, you won. After this we can start talking about how to write code which is of a high standard.

For PhD students and more mathematically mature people who are interested in seeing what's going on, it's always worthwhile pointing out current mathematical projects which are run via Zulip at {\tt leanprover.zulipchat.com}, the platform behind the Lean community chat. We are always looking for new people to help out with various projects, and are happy to ``onboard'' newcomers. Right now the three main active projects on the site are a formalization of the proof of Fermat's Last Theorem for regular primes being led by Riccardo Brasca, a formalization of sphere eversion being led by Patrick Massot, and ``the liquid tensor experiment'', being led by Johan Commelin, with a lot of help from Adam Topaz. The latter is the formalization of a 2019 theorem of Clausen and Scholze; Scholze challenged the formalization community to verify a crucial theorem he had announced with Clausen, and the Lean community rose to the challenge. The back story is interesting; Scholze suggested that perhaps the current refereeing process that we have in mathematics would not dive deeply into a specific part of the work, and Scholze was interested to see whether a theorem prover could do this instead. It turns out that this was indeed possible. Scholze was an advisor throughout the project; the main technical intermediate result is completely formalized, and the actual challenge itself is very nearly completed.

\subsection*{Teaching formalization skills}

Teaching students how to formalize mathematics means teaching them how to translate mathematical ideas between English and Lean. Just like learning a foreign human language, you begin by translating basic stuff between the two languages, and you ask if you don't understand something. In the Lean course I teach in Imperial's mathematics department, we work through human mathematics that the students have already been taught, and we learn to speak Lean's language. Later on I introduce some new mathematics, when the students have already learnt some Lean. Athina Thoma and Paola Iannone taught me that teaching the first years equivalence relations and Lean at the same time would usually not end well. However to a student who knows both, formalising the fact that equivalence relations biject with partitions would be an excellent Lean learning experience. The art is to find the right project, and the correct project typically depends on the student.

\subsection*{Fixing issues in modern research}

In an interview with Wired I was once quoted as suggesting that all maths was wrong. This is something I know not to be true -- some maths is definitely correct, and (at least in classical logic) these statements are opposite to one another. However, I did say that perhaps some of our castles were built on sand, and I now think even this is a bit naive. Having had conversations with David Rabouin, a historian and philosopher of mathematics, I now understand that cutting edge maths always looks like this. There are bits which are not quite checked but everyone knows it will usually work out fine, or at least well enough to make the main theorem go through.

Sometimes maths does go wrong though. This century, in my area alone (number theory), prestigious mathematicians have announced proofs of Leopoldt's conjecture and the ABC conjecture; the latter work was even published in a reputable mathematical journal. Yet our community does not seem to accept the proofs as rigorous. Unfortunately, Lean will not deliver us from these problems, at least not yet. Right now, Lean proves theorems by having humans translate those theorems from the English, and if nobody is prepared to take on the monumental task of translating the Mochuzuki proof of the ABC conjecture into a theorem prover (and why should anybody? This is not how the mathematical community has treated published proofs in the past), then I don't see any way past the impasse.

Some proofs have already gone beyond the one brain barrier -- no one human understands all of the ideas in it. Areas which were formally fashionable can die out if the big conjectures driving them forwards are resolved. Anything not properly documented actually runs the risk of being lost. One can hope that new and simpler proofs will come along. But this is not always true. Whether or not we choose to use theorem provers to do it, mathematicians need to start thinking more carefully than ever about precisely documenting what we think we already know, so we can answer technical questions from future generations of mathematicians.

%% file: 8-paulson.tex
\section{When will computer-assisted proof become part of everyday mathematics?}\label{sec:paulson}
\secauthor{Lawrence Paulson}

\vspace{-\baselineskip}
\subsection*{Introduction and Background}
The idea of applying technology to mathematical reasoning began to be realized in the 1960s. 
N.~G. de~Bruijn's AUTOMATH~\cite{automath} was a type theory for expressing mathematical definitions and proofs.  Trybulec's Mizar system~\cite{Grabowski2015} included a human-readable formal language for abstract mathematics.  
Both were professional mathematicians and intended their work to be beneficial to mathematics itself, but the technology was not ready. Nevertheless, they made valuable progress. AUTOMATH led to today's dependent type theories. Mizar accumulated a substantial and wide-ranging library of formalized mathematics, while its readable structured language is still the best there is.

\textit{Interactive theorem provers}, or proof assistants, emerged during the 1980s from LCF and HOL through the work of Milner, Gordon, et al.~\cite{mgordon-history,paulson-computational}. These tools were intended for verification in computer science. HOL (for higher-order logic) had been chosen in order to implement a certain style of hardware verification. These verification tasks seldom required any mathematics beyond the integers.

The 1994 discovery of a bug in the floating-point unit of Intel's Pentium dramatically focused the verification community on the real numbers. 
John Harrison formally proved the correctness of an algorithm for computing the exponential function, taking account of all the peculiarities of floating-point arithmetic~\cite{harrison-exp}. He went on to play a major role in the \textit{Flyspeck project}: the formal confirmation of Thomas Hales's proof of the Kepler conjecture~\cite{hales-formal-Kepler}. He formalized many landmark results in mathematics, such as the prime number theorem~\cite{harrison-pnt}.

Gonthier's formalization of the Four Color Theorem \cite{gonthier-4ct} had already demonstrated that an interactive theorem prover (Coq in this case) could help settle a genuine question in mathematics. This task was similar to Flyspeck in that it involved formally checking a large number of computations.  

\subsection*{Isabelle}

Isabelle \cite{paulson-isa-book} emerged from the LCF tradition with the aim of supporting a multiplicity of formalisms, including set theory. However, with higher-order logic dominating the verification world, Isabelle/HOL became the dominant instance of Isabelle. Its formalism extended that of the various HOLs with axiomatic type classes, allowing the systematic reuse of formal material sharing the same axiomatic basis~\cite{hoelzl-filters}. Isabelle adopted a structured proof language, Isar, based on Mizar's mathematical language~\cite{wenzel-isabelle/isar}. Isar expresses proofs with a nested structure where milestones---intermediate claims and proofs---are explicitly written out. 
Isabelle provides automation for proof (powerful external provers can be invoked through \textit{sledgehammer})
and also for disproof in the form of \textit{counterexample search}~\cite{blanchette-nitpick}. The user interface is a unique Interactive Development Environment for editing live proof documents~\cite{wenzel-pide}. 

There are powerful synergies between these features. Structured proof text is easy to reuse. Copied into a new development, it will instantly be checked and any errors flagged. Structured proofs also work well with sledgehammer: if a given statement is too difficult to be proved automatically, the user may propose an intermediate statement that might be easy enough to prove automatically and lead eventually to a proof of the original statement.

This powerful automation has greatly relieved the tedium that accompanied formal proof in the early days. It has also helped encourage the growth of substantial libraries of formalized mathematics. 

\subsection*{Proof Assistants: Ready for Prime Time?}

Today we have a wide variety of systems: the HOL family, for higher-order logic; Isabelle/HOL, for higher-order logic with axiomatic type classes; Coq and Lean, implementing a dependent type theory. These choices are motivated by conflicting priorities such as the simplicity of the implementation and its semantics as opposed to the expressiveness of the logical calculus.
Competition between research groups was spurred by Freek Wiedijk's online list detailing which of the ``top 100 theorems''\footnote{\url{http://www.cs.ru.nl/~freek/100/}} had been proved in various systems. Only three of the hundred theorems remain unproved, and one of those is Fermat's Last Theorem!

Much of modern mathematics falls within the scope of existing verification tools. 
There are now huge libraries of formalized mathematics. Isabelle's \textit{Archive of Formal Proofs}\footnote{\url{http://www.isa-afp.org}} 
contains over 650 entries including includes much core mathematics---linear algebra, multivariate analysis, probability, complex analysis, topology---and over 3 million lines of proofs.
Lean's \textit{mathlib} is an immense and rapidly-growing corpus of material from every branch of mathematics.%
\footnote{\url{https://leanprover-community.github.io/mathlib_docs/}}

So are these systems finally ready to support mathematicians? There are still many obstacles:
\begin{itemize}
  \item Formal syntax looks artificial and often is barely legible. To see the difficulties presented by traditional notation, contrast $x^2$, $\nabla^2f$, $\sin^2 \theta$, $f^2 (x)$. In group theory, $G$ denotes a group but also a set; $ab$ is the product of $a$ and $b$ but $HK$ is quite a different thing, simply because we used a different part of the alphabet. 
In set theory, $\lambda < \aleph_1$ has a different meaning from $\lambda < \omega_1$ even though $\aleph_1=\omega_1$.
  \item The libraries of formalized mathematics still have many gaps, and what they do have is difficult to find. The names of theorems are frequently ambiguous (what is Roth's theorem?), as are concepts such as limits (which could refer to analysis, topology or even set theory). But one should be able to search for ``limit theorems''  and get something relevant. Better still, the system might make suggestions unprompted.
  \item Obvious statements are often too hard to prove. Examples include showing a set to be finite, showing a function to be continuous, calculating a derivative and evaluating a limit. In some cases, skilful use of the existing automation can be effective. In others, specialized decision procedures must be coded. Best known are decision procedures for linear arithmetic; a recent development is Eberl's limit solver~\cite{eberl-verified-real}. 
  \end{itemize}
  
\subsection*{Toward the Future}
  
Recently, with funding from the European Research Council,\footnote{Project ALEXANDRIA, GA 742178} my colleagues and I have been exploring and stretching the limits of today's technology. We have formalized relatively recent (post 1970) and deep results into both classical~\cite{koutsoukou-irrationality} and previously overlooked fields. The latter include combinatorial design theory~\cite{edmonds-modular-first} and ordinal partition theory~\cite{dzamonja-formalising}. We have formalized some sophisticated definitions, notably Grothendieck schemes---necessary for advanced work in algebraic geometry---which had hitherto been thought to lie outside the scope of Isabelle/HOL's simple type theory~\cite{bordg-simple-tr}.  

Mathematicians also need help navigating our huge library of formalized mathematics. It can be hard to know whether a desired result has been formalized: many theorems go by various names, or conversely, one name (e.g.\ Young's inequality) may be applied to a family of distinct results. For every well-known result, there may be dozens of technical lemmas, mostly obvious and yet needed in order to prove anything on today's systems. My colleagues have been building an experimental search engine, called SErAPIS,\footnote{\url{https://behemoth.cl.cam.ac.uk/search/}} that makes it possible to search the entire library with the help of a huge dictionary of mathematical concepts.

Also attractive is the idea that the proofs in our library might be used to generate new proofs automatically. This is another way the computer can help the user get value out of our 3 million lines of formal proofs. Promising results are just starting to appear from a number of research groups.

For the further future, we may hope to mechanize mathematical intuition. This is the knowledge that tells us that a given function is surely continuous or that a certain formula cannot generate only prime numbers. That's how we know that a particular claim cannot be true as stated or proved using the methods advertised. Formalization has demonstrated time and again that while published proofs often contain errors, theorem statements are generally correct. Mathematicians can perceive the truth even if they can't always write down the correct argument. Giving such intuition to a computer would transform our field. This task will remain open for the next generation or two.

%% file: main.bbl
\begin{thebibliography}{10}

\bibitem{1977SciAm.237d.108A}
K.~{Appel} and W.~{Haken}.
\newblock {The Solution of the Four-Color-Map Problem}.
\newblock {\em Scientific American}, 237:108--121, Oct. 1977.

\bibitem{Baker70}
I.~N. Baker.
\newblock Completely invariant domains of entire functions.
\newblock {\em Mathematical Essays Dedicated to A. J. Macintyre.}, pages
  33--35, 1970.

\bibitem{bayer-beginners-quest}
J.~Bayer, M.~David, A.~Pal, and B.~Stock.
\newblock Beginners' quest to formalize mathematics: A feasibility study in
  isabelle.
\newblock In C.~Kaliszyk, E.~Brady, A.~Kohlhase, and C.~Sacerdoti~Coen,
  editors, {\em Intelligent Computer Mathematics}, pages 16--27, Cham, 2019.
  Springer International Publishing.

\bibitem{Digit_Expansions-AFP}
J.~Bayer, M.~David, A.~Pal, and B.~Stock.
\newblock Digit expansions.
\newblock {\em Archive of Formal Proofs}, Apr. 2022.
\newblock \url{https://isa-afp.org/entries/Digit_Expansions.html}, Formal proof
  development.

\bibitem{dprm_isabelle}
J.~Bayer, M.~David, A.~Pal, B.~Stock, and D.~Schleicher.
\newblock {The DPRM Theorem in Isabelle (Short Paper)}.
\newblock In J.~Harrison, J.~O'Leary, and A.~Tolmach, editors, {\em 10th
  International Conference on Interactive Theorem Proving (ITP 2019)}, volume
  141 of {\em Leibniz International Proceedings in Informatics (LIPIcs)}, pages
  33:1--33:7, Dagstuhl, Germany, 2019. Schloss Dagstuhl--Leibniz-Zentrum fuer
  Informatik.

\bibitem{DPRM_Theorem-AFP}
J.~Bayer, M.~David, B.~Stock, A.~Pal, Y.~Matiyasevich, and D.~Schleicher.
\newblock Diophantine equations and the {DPRM} theorem.
\newblock {\em Archive of Formal Proofs}, June 2022.
\newblock \url{https://isa-afp.org/entries/DPRM_Theorem.html}, Formal proof
  development.

\bibitem{J41}
C.~Benzm{\"u}ller.
\newblock Universal (meta-)logical reasoning: Recent successes.
\newblock {\em Science of Computer Programming}, 172:48--62, March 2019.
\newblock Preprint: \url{http://dx.doi.org/10.13140/RG.2.2.11039.61609/2}.

\bibitem{J40}
C.~Benzm{\"u}ller and D.~S. Scott.
\newblock Automating free logic in {HOL}, with an experimental application in
  category theory.
\newblock {\em Journal of Automated Reasoning}, 2019.
\newblock Preprint: \url{http://dx.doi.org/10.13140/RG.2.2.11432.83202}.

\bibitem{C55}
C.~Benzm{\"u}ller and B.~Woltzenlogel~Paleo.
\newblock The inconsistency in {G{\"o}del's} ontological argument: A success
  story for {AI} in metaphysics.
\newblock In S.~Kambhampati, editor, {\em IJCAI 2016}, volume 1-3, pages
  936--942. AAAI Press, 2016.

\bibitem{Coq}
Y.~Bertot and P.~Cast{\'e}ran.
\newblock {\em Interactive theorem proving and program development. Coq’Art:
  the calculus of inductive constructions}.
\newblock Springer-Verlag Berlin Heidelberg, 2013.

\bibitem{blanchette-nitpick}
J.~C. Blanchette and T.~Nipkow.
\newblock {Nitpick}: A counterexample generator for higher-order logic based on
  a relational model finder.
\newblock In M.~Kaufmann and L.~C. Paulson, editors, {\em Interactive Theorem
  Proving}, volume 6172 of {\em Lecture Notes in Computer Science}, pages
  131--146. Springer, 2010.

\bibitem{bordg-simple-tr}
A.~Bordg, L.~Paulson, and W.~Li.
\newblock Simple type theory is not too simple: {Grothendieck's} schemes
  without dependent types.
\newblock {\em Experimental Mathematics}, 2022.

\bibitem{Lean10}
M.~Carneiro.
\newblock A {L}ean formalization of {M}atiyasevič's theorem, 2018.
\newblock arxiv: \url{https://arxiv.org/abs/1802.01795}.

\bibitem{liquid_tensor_nature}
D.~Castelvecchi.
\newblock Mathematicians welcome computer-assisted proof in ‘grand
  unification’theory.
\newblock {\em Nature}, 595(7865):18--19, 2021.

\bibitem{dzamonja-formalising}
M.~D{\v z}amonja, A.~Koutsoukou-Argyraki, and L.~C. Paulson.
\newblock Formalising ordinal partition relations using {Isabelle/HOL}.
\newblock {\em Experimental Mathematics}, 2021.

\bibitem{eberl-verified-real}
M.~Eberl.
\newblock Verified real asymptotics in {Isabelle/HOL}.
\newblock In {\em International Symposium on Symbolic and Algebraic Computation
  (ISAAC 2019)}, pages 147--154, 2019.

\bibitem{edmonds-modular-first}
C.~Edmonds and L.~C. Paulson.
\newblock A modular first formalisation of combinatorial design theory.
\newblock In F.~Kamareddine and C.~Sacerdoti~Coen, editors, {\em Intelligent
  Computer Mathematics}, pages 3--18, Cham, 2021. Springer International
  Publishing.

\bibitem{Gonthier08}
G.~Gonthier.
\newblock {Formal Proof -- The Four-Color Theorem}.
\newblock {\em Notices of the American Mathematical Society},
  55(11):1382--1393, Dec. 2008.

\bibitem{gonthier-4ct}
G.~Gonthier.
\newblock The four colour theorem: Engineering of a formal proof.
\newblock In D.~Kapur, editor, {\em Computer Mathematics}, LNCS 5081, pages
  333--333. Springer, 2008.

\bibitem{mgordon-history}
M.~J.~C. Gordon.
\newblock From {LCF} to {HOL}: a short history.
\newblock In G.~Plotkin, C.~Stirling, and M.~Tofte, editors, {\em Proof,
  Language, and Interaction: Essays in Honor of {Robin Milner}}, pages
  169--185. MIT Press, 2000.

\bibitem{Grabowski2015}
A.~Grabowski, A.~Korni{\l}owicz, and A.~Naumowicz.
\newblock Four decades of {Mizar}.
\newblock {\em Journal of Automated Reasoning}, 55(3):191--198, Oct 2015.

\bibitem{hales-formal-Kepler}
T.~Hales, M.~Adams, G.~Bauer, T.~D. Dang, J.~Harrison, L.~T. Hoang,
  C.~Kaliszyk, V.~Magron, S.~Mclaughlin, T.~T. Nguyen, et~al.
\newblock A formal proof of the {Kepler} conjecture.
\newblock {\em Forum of Mathematics, Pi}, 5:e2, 2017.

\bibitem{harrison-exp}
J.~Harrison.
\newblock Floating point verification in {HOL} {L}ight: the exponential
  function.
\newblock {\em Formal Methods in System Design}, 16:271--305, 2000.

\bibitem{harrison-pnt}
J.~Harrison.
\newblock Formalizing an analytic proof of the prime number theorem.
\newblock {\em Journal of Automated Reasoning}, 43(3):243--261, 2009.

\bibitem{hartnett2021}
K.~Hartnett.
\newblock Proof assistant makes jump to big-league math.
\newblock Online at
  \url{https://www.quantamagazine.org/lean-computer-program-confirms-peter-scholze-proof-20210728/},
  July 2021.

\bibitem{Heule2018}
M.~J.~H. Heule.
\newblock Schur number five.
\newblock In {\em {AAAI}}, pages 6598--6606. {AAAI} Press, 2018.

\bibitem{Heule2017}
M.~J.~H. Heule and O.~Kullmann.
\newblock The science of brute force.
\newblock {\em Commun. ACM}, 60(8):70--79, July 2017.

\bibitem{hoelzl-filters}
J.~H{\"o}lzl, F.~Immler, and B.~Huffman.
\newblock Type classes and filters for mathematical analysis in {Isabelle/HOL}.
\newblock In S.~Blazy, C.~Paulin-Mohring, and D.~Pichardie, editors, {\em
  Interactive Theorem Proving}, LNCS 7998, pages 279--294. Springer, 2013.

\bibitem{digmath}
{IMU Committee on Electronic Information and Communication ({CEIC})}.
\newblock World digital mathematics library ({WDML}).
\newblock
  \url{https://www.mathunion.org/ceic/library/world-digital-mathematics-library-wdml}.

\bibitem{J50}
D.~Kirchner, C.~Benzm{\"u}ller, and E.~N. Zalta.
\newblock Mechanizing principia logico-metaphysica in functional type theory.
\newblock {\em Review of Symbolic Logic}, 13(1):206--218, 2020.

\bibitem{koutsoukou-irrationality}
A.~Koutsoukou-Argyraki, W.~Li, and L.~C. Paulson.
\newblock Irrationality and transcendence criteria for infinite series in
  {Isabelle/HOL}.
\newblock {\em Experimental Mathematics}, 2021.

\bibitem{Lamport_ErrorsInProofs}
L.~Lamport.
\newblock Errors in proofs — a correction and further data.
\newblock \url{https://lamport.azurewebsites.net/tla/proof-statistics.html}.

\bibitem{Lamport_HowToWriteProof}
L.~Lamport.
\newblock How to write a proof.
\newblock {\em The American mathematical monthly}, 102(7):600--608, 1995.
\newblock Also appeared in Karen Uhlenbeck, editor, \emph{Global Analysis in
  Modern Mathematics}, Houston, 1992. Publish or Perish Press.

\bibitem{Lamport_HowToWrite21}
L.~Lamport.
\newblock How to write a 21 st century proof.
\newblock {\em Journal of fixed point theory and applications}, 11(1):43--63,
  2012.

\bibitem{form3}
D.~Larchey-Wendling and Y.~Forster.
\newblock {Hilbert's Tenth Problem in Coq (Extended Version)}.
\newblock {\em {Logical Methods in Computer Science}}, {Volume 18, Issue 1},
  Mar. 2022.

\bibitem{Lasse_Dave}
D.~S. Lasse Rempe-Gillen.
\newblock On connected preimages of simply-connected domains under entire
  functions.
\newblock 2018.
\newblock Manuscript. arxiv: \url{https://arxiv.org/pdf/1801.06359.pdf}.

\bibitem{matobvo2}
Y.~V. Matiyasevich.
\newblock Mathematical proof: yesterday, today, tomorrow.
\newblock Joint meeting of the St. Petersburg Mathematical Society and the
  Section of the House of Scientists.
  \url{http://www.mathnet.ru/php/seminars.phtml?presentid=2036&option_lang=eng}.

\bibitem{automath}
R.~P. Nederpelt, J.~H. Geuvers, and R.~C. de~Vrijer, editors.
\newblock {\em Selected Papers on Automath}.
\newblock North-Holland, 1994.

\bibitem{form1a}
K.~Pak.
\newblock The matiyasevich theorem. preliminaries.
\newblock 25(4):315--322.

\bibitem{paulson-isa-book}
L.~C. Paulson.
\newblock {\em Isabelle: A Generic Theorem Prover}.
\newblock Springer, 1994.
\newblock LNCS 828.

\bibitem{paulson-computational}
L.~C. Paulson.
\newblock Computational logic: Its origins and applications.
\newblock {\em Proceedings of the Royal Society of London A: Mathematical,
  Physical and Engineering Sciences}, 474(2210), 2018.

\bibitem{fm}
{R. Matuszewski (editor in-chief)}.
\newblock Formalized mathematics.
\newblock \url{https://fm.mizar.org/}.

\bibitem{Spivak}
M.~Spivak.
\newblock {\em Calculus}.
\newblock W. A. Benjamin, Inc., New York, 1967.

\bibitem{Thurston94}
W.~P. Thurston.
\newblock On proof and progress in mathematics.
\newblock {\em Bulletin of the American Mathematical Society}, 30(2):161–177,
  Apr 1994.

\bibitem{Vavilov}
N.~Vavilov.
\newblock Reshaping the metaphor of proof.
\newblock {\em Philosophical Transactions of the Royal Society A},
  377(2140):20180279, 2019.

\bibitem{matobvo1}
N.~A. Vavilov.
\newblock Mathematical proof: yesterday, today, tomorrow.
\newblock Joint meeting of the St. Petersburg Mathematical Society and the
  Section of the House of Scientists.
  \url{http://www.mathnet.ru/php/seminars.phtml?presentid=2035&option_lang=eng}.

\bibitem{matobvo3}
A.~M. Vershik.
\newblock Mathematical proof: yesterday, today, tomorrow.
\newblock Joint meeting of the St. Petersburg Mathematical Society and the
  Section of the House of Scientists.
  \url{http://www.mathnet.ru/php/seminars.phtml?presentid=2583&option_lang=eng}.

\bibitem{voevodsky-origins}
V.~Voevodsky.
\newblock The origins and motivations of univalent foundations.
\newblock {\em The Institute Letter}, pages 8--9, Summer 2014.
\newblock Online at \url{https://www.ias.edu/ideas/2014/voevodsky-origins}.

\bibitem{wenzel-isabelle/isar}
M.~Wenzel.
\newblock {Isabelle/Isar} --- a generic framework for human-readable proof
  documents.
\newblock {\em Studies in Logic, Grammar, and Rhetoric}, 10(23):277--297, 2007.
\newblock From Insight to Proof --- Festschrift in Honour of Andrzej Trybulec.

\bibitem{wenzel-pide}
M.~Wenzel.
\newblock Asynchronous user interaction and tool integration in
  {Isabelle/PIDE}.
\newblock In G.~Klein and R.~Gamboa, editors, {\em Interactive Theorem Proving
  --- 5th International Conference, {ITP} 2014}, LNCS 8558, pages 515--530.
  Springer, 2014.

\bibitem{100}
F.~Wiedijk.
\newblock Formalizing 100 theorems.
\newblock \url{http://www.cs.ru.nl/~freek/100/}.

\bibitem{compinstr}
Ю.В.Матиясевич.
\newblock Математическое доказательство:
  вчера, %сегодня, завтра.
\newblock \emph{Компьютернве инструменты в
  образовании}, выпуск 6, cтр. 13--24, 2012 г. This is a
  printed report of the talk \cite{matobvo2}.

\end{thebibliography}
